\documentclass[leqno,12pt]{article}
\usepackage{latexsym}
\usepackage{amsmath}
\usepackage{amssymb}
\usepackage{bbm}
\usepackage{color}
\usepackage[numbers]{natbib}

\usepackage{graphicx}
\usepackage{epsfig}

 0

\newtheorem{thm}{Theorem}[section]

\newcommand{\ea}{\end{array}}
\newcommand{\beqohne}{\begin{eqnarray*}}
\newcommand{\eeqohne}{\end{eqnarray*}}
\newcommand{\beohne}{\begin{equation*}}
\newcommand{\eeohne}{\end{equation*}}

\def\proof{\noindent{\bf Proof:}\hskip10pt}

\def \sur#1#2{\mathrel{\mathop{\kern 0pt#1}\limits^{#2}}}

\makeatletter\@addtoreset{equation}{section}\makeatother

\def\proof{\noindent{\bf Proof:}\hskip10pt}

\begin{document}

\title{Addendum to "Sum rules via large deviations"}


\author{{\small Fabrice Gamboa}\footnote{ Universit\'e Paul Sabatier, Institut de Math\'ematiques de Toulouse,  31062 Toulouse Cedex 9, France, 
e-mail: gamboa@math.univ-toulouse.fr}
\and{\small Jan Nagel}\footnote{Technische Universitat M\"unchen, Fakult\"at f\"ur Mathematik, Boltzmannstr. 3, 85748 Garching, Germany,  e-mail: jan.nagel@tum.de}
\and{\small Alain Rouault}\footnote
{Laboratoire de Math\'ematiques de Versailles, UVSQ, CNRS, Université Paris-Saclay, 78035-Versailles Cedex France, e-mail: alain.rouault@uvsq.fr}}

\maketitle

\abstract{In these notes we fill a gap in a proof in Section 4 of Gamboa, Nagel, Rouault [Sum rules via large deviations,  J. Funct. Anal. 270 (2016), 509-559]. 
We prove a general theorem which combines a LDP with a convex rate function and a LDP with a non-convex one. This result 
will be used to prove LDPs for spectral matrix measures and for spectral measures on the unit circle. 
}
\bigskip

\smallskip


\section{Introduction}

In Section 4 of \cite{magicrules}, we studied  large deviations for a pair of random variables with values in topological vector spaces by means of the joint normalized generating function. However, in some cases, the rate function of one of the marginals is not convex, which invalidates this way of proof.  
Actually, it is possible to state  a general theorem which combines a LDP with a convex rate function and a LDP with a non-convex one. It is used to prove LDPs in \cite{GaNaRomat} for spectral matrix measures  and in \cite{GNROPUC} for spectral measures on the unit circle.\\
Since this theorem may have its own interest, we give it in a general setting in Section \ref{general} after recalling classical results in Section \ref{class}. We 
come back in Section \ref{appli} to the framework of  \cite{magicrules}.

In the sequel, we assume that $\mathcal X$ and $\mathcal Y$ are Hausdorff topological vector spaces. $\mathcal X^*$ is the topological dual of $\mathcal X$ and $\mathcal{X}$ is endowed with the weak topology. We denote by $C_b(\mathcal Y )$ the set of all bounded continuous functions $\varphi:\mathcal Y \to \mathbb{R}$. 
A point $x\in \mathcal{X}$ is called an exposed point of a function $F$ on $\mathcal{X}$, if there exists $x^*\in \mathcal{X}^*$ (called an exposing hyperplane for $x$) such that 
\begin{align}\label{exposinghyper}
F(x) - \langle x^*,x\rangle < 
F(z) - \langle x^*,z\rangle
\end{align}
for all $z\neq x$.

\section{Some classical results in large deviations}
\label{class}
Let us recall two well known results in the theory of large deviations which will be combined  in order to solve our problem. The first result is the inverse of Varadhan's lemma (Theorem 4.4.2 in \cite{demboz98}), the second one is a version of the so-called Baldi's theorem  (Theorem 4.5.20 in \cite{demboz98}). The latter differs from the version in \cite{demboz98} in a straightforward condition to identify the rate function, which was applied for instance 
 in \cite{grz} (see also \cite{gamb}).
The proof of our Theorem \ref{newgeneral} will be quite similar to the 
proof of these two classical theorems.

\begin{thm}[Bryc's Inverse Varadhan Lemma]
\label{Bryc}
Suppose that the sequence  $(Y_n)$ of random variables  in $\mathcal{Y}$ is exponentially tight and that the limit
\[\Lambda (\varphi) := \lim_{n\to \infty} \frac{1}{n} \log \mathbb E e^{n\varphi(Y_n)}\] exists for every $\varphi \in \mathcal C_b(\mathcal Y)$. Then $(Y_n)$ satisfies the LDP with the good rate function 
\[\mathcal I(y) = \sup_{\varphi \in \mathcal C_b(\mathcal{Y})} \{\varphi(y) - \Lambda(\varphi)\}\,.\]
Furthermore, for every $\varphi \in \mathcal C_b(\mathcal{Y})$,
\[\Lambda(\varphi) = \sup_{y\in \mathcal Y} \{\varphi(y) - \mathcal I(y)\}\,.\]
\end{thm}

\begin{thm}[A version of Baldi's Theorem]
\label{Baldi}
Suppose that the sequence  $(X_n)$ of random variables in $\mathcal X$ is exponentially tight and that   : 
\begin{enumerate}
\item There is a  set $D\subset \mathcal X^*$ and a function $G_X:D\to \mathbb{R}$ such that for all $x^*\in D$ 
\begin{equation}
\label{ncgf}\lim_{n\to \infty} \frac{1}{n} \log \mathbb E \exp \left(n \langle x^* , X_n\rangle\right) =  G_X(x^*) \,.\end{equation}
\item If $\mathcal{F}$ denotes the set of exposed points $x$ of
\[G_X^* (x) = \sup_{x^* \in D} \{\langle x^*,x\rangle - G_X(x^*)\}  , \]
with an exposing hyperplane $x^*$ satisfying $x^*\in D$ and 
$\gamma x^*\in D$ for some $\gamma
>1 $, then for every $x\in \{ G_X^* < \infty\}$ there exists a sequence $(x_k)_k$ with $x_k \in \mathcal{F}$ such that $\lim_{k\to \infty} x_k = x$ and 
\begin{align*}
\lim_{k\to \infty} G_X^*(x_k) =G_X^*(x) .
\end{align*}
\end{enumerate}
Then $(X_n)$ satisfies the LDP with good rate function $G_X^*$.  
\end{thm}

\section{A general theorem}
\label{general}
Our extension is the following combination of the two above theorems. The main point is that the rate function does not need to be convex, but we still only need to control linear functionals of $X_n$. 
\begin{thm}
\label{newgeneral}
Assume that $X_n \in \mathcal X$ and $Y_n \in \mathcal Y$ are defined on the same probabilistic space. Moreover, we assume that the two sequences $(X_n)$ and $(Y_n)$ are exponentially tight. Assume further that:
\begin{enumerate}
\item There is a set $D\subset \mathcal X^*$ and functions $G_X:D\to \mathbb{R}$, $J:C_b(\mathcal Y) \to \mathbb{R}$ such that
for all $x^*\in D$ and $\varphi\in C_b(\mathcal Y)$
\begin{equation}
\label{newncgf}\lim_{n\to \infty} \frac{1}{n} \log \mathbb E \exp \left(n \langle x^* , X_n\rangle + n \varphi(Y_n)\right) =  G_X(x^*) + J (\varphi)\,.
\end{equation}
\item
If $\mathcal{F}$ denotes the set of exposed points $x$ of
\[G_X^* (x) = \sup_{x^* \in D} \{\langle x^*,x\rangle - G_X(x^*)\}   \]
with an exposing hyperplane $x^*$ satisfying $x^*\in D$ and 
$\gamma x^*\in D$ for some $\gamma
>1 $, then for every $x\in \{ G_X^* < \infty\}$ there exists a sequence $(x_k)_k$ with $x_k \in \mathcal{F}$ such that $\lim_{k\to \infty} x_k = x$ and 
\begin{align*}
\lim_{k\to \infty} G_X^*(x_k) =G_X^*(x) .
\end{align*}
\end{enumerate}  
Then, the pair $(X_n, Y_n)$ satisfies the LDP with speed $n$ and good rate function
\[\mathcal I(x,y) = G_X^*(x) + \mathcal{I}_Y(y)\,,\]
where 
\begin{equation*}
\mathcal{I}_Y(y) = \sup_{\varphi \in C_b(\mathcal Y)} \{ \varphi(y) - J(\varphi)\} .
\end{equation*}
\end{thm}
Let us note that in view of Varadhan's Lemma we have
\[J(\varphi) = \sup_{y \in \mathcal Y} \{ \varphi(y) - \mathcal{I}_Y(y) \}. \]

\proof

{\bf Upperbound}: The proof  follows the lines of the proof of part (b) of Theorem 4.5.3 in \cite{demboz98}. Note that since the sequence $(X_n,Y_n)$ is exponentially tight it suffices to show the upper bound for compact sets. Furthermore, the rate is necessarily good, since, if in \eqref{newncgf} we set $x^*=0$ (resp. $\varphi=0$) it reduces to Theorem \ref{Bryc} (resp. Theorem \ref{Baldi}).

{\bf Lowerbound}: As usual, it is enough to consider a neighbourhood 
 $\Delta_1 \times \Delta_2$ of $(x,y)$ where  $\mathcal I(x,y) < \infty$, take 
$\liminf_{n \to \infty} \frac{1}{n} \log \mathbb P ((X_n , Y_n) \in \Delta_1 \times \Delta_2)$ and get a lower bound tending to $\mathcal I (x,y)$  when the infimum over all neighborhoods is taken. Actually, due to the density assumption 2. it is enough to study the  lower bound of $\mathbb P(X_n \in \Delta_1, Y_n \in \Delta_2)$  when 
 $x\in \mathcal F$ and $\mathcal I_Y (y) < \infty$. 

As in \cite{demboz98} (Proof of Lemma 4.4.6), let $\varphi:\mathcal Y \rightarrow [0,1]$ be a continuous function, such that $\varphi (y) =1$ and $\varphi$ vanishes on the complement $\Delta_2^c$ of $\Delta_2$. For $m >0$, define $\varphi_m := m(\varphi - 1)$. Note that
\[J(\varphi_m) \geq -\mathcal{I}_Y(y)\,.\]
We have 
\begin{eqnarray*}
\mathbb P(X_n \in \Delta_1, Y_n \in \Delta_2) =
\mathbb E \left[ {\mathbbm{1} }_{\{ X_n \in \Delta_1\}} {\mathbbm{1}}_{\{ Y_n \in \Delta_2\}}e^{ n\langle x^*, X_n\rangle+n\varphi_m(Y_n)} e^{-n\langle x^*, X_n\rangle-n\varphi_m(Y_n)}\right] .
\end{eqnarray*}
Now $-\varphi_m \geq 0$ and on $\Delta_1$, $-\langle x^*, X_n\rangle \geq -\langle x^*, x\rangle -  \delta$ for a $\delta>0$, so that
\begin{equation}
\mathbb P(X_n \in \Delta_1, Y_n \in \Delta_2) \geq \mathbb E  \left[{\mathbbm{1} }_{\{ X_n \in \Delta_1\}} {\mathbbm{1}}_{\{ Y_n \in \Delta_2\}} e^{ n\langle x^*, X_n\rangle+n\varphi_m(Y_n)} \right] e^{- n\langle x^*, x\rangle - n \delta}\,.
\end{equation} 
Denoting
\[\ell_n = \frac{1}{n}  \log \mathbb E  e^{ n\langle x^*, X_n\rangle} \ , \ \mathcal L_n := \frac{1}{n} \log \mathbb E e^{n\langle x^*, X_n\rangle + n \varphi_m (Y_n) }\]
and 
$\widetilde{\mathbb P}$ the new probability on $\mathcal X \times \mathcal Y$ such that
\[\frac{d\widetilde{\mathbb P}}{d\mathbb P}=  e^{n\langle x^*, X_n\rangle + n \varphi_m (Y_n)  -n \mathcal L_n}\,,\]
we get
\begin{equation}
\label{mino}
\mathbb P(X_n \in \Delta_1, Y_n \in \Delta_2) \geq \widetilde{\mathbb P}(X_n \in \Delta_1, Y_n \in \Delta_2)e^{- n\langle x^*, x\rangle - n \delta + n \mathcal L_n}\,.
\end{equation}
For the exponential term we have
\begin{equation}
\label{uptodelta}
\liminf_{n\to \infty} \frac{1}{n} \log e^{- n\langle x^*, x\rangle - n \delta + n \mathcal L_n}\geq \langle x^*, x\rangle -  \delta + G_X
(x^*) + J(\varphi_m) \geq -
G_X^*(x)  -\mathcal{I}_Y(y)-\delta .
\end{equation}
We may choose $\delta $ arbitrarily small by choosing $\Delta_1$ sufficiently small, 
so that it will be enough to prove that
\begin{equation}
\label{7.7}
\widetilde{\mathbb P}(X_n \in \Delta_1, Y_n \in \Delta_2) \xrightarrow[n\to \infty]{} 1
\end{equation}
or equivalently, that
\begin{equation} \label{complementprob}
\widetilde{\mathbb P}(X_n \in \Delta_1^c) + \widetilde{\mathbb P}( Y_n \in \Delta_2^c) \xrightarrow[n\to \infty]{} 0\,.
\end{equation}
For the first term, note that under $\widetilde{\mathbb P}$ the moment generating function of $X_n$ satisfies
\begin{align*}
\lim_{n\to \infty} \frac{1}{n} \log \widetilde{\mathbb{E}} [e^{n\langle z^*,X_n\rangle}] 
& = \lim_{n\to \infty} \frac{1}{n} \log {\mathbb{E}} [e^{n\langle z^*+x^*,X_n\rangle + \varphi_m(Y_n)-n\mathcal{L}_n}] \\
& = G_X(z^*+x^*) + J(\varphi_m) - G_X(x^*) - J(\varphi_m) \\
& = G_X(z^*+x^*) - G_X(x^*)\\
& = : \widetilde G_X(z^*) ,
\end{align*}
for $z^* \in \widetilde D  := \{z^* : x^* + z^* \in D \}$. 
We may then  follow the argument on  p.159-160 in \cite{demboz98} (as an auxiliary result in  their proof of the lower bound).
Using that $x^*\in D$ is an exposing hyperplane, we get 
\begin{align*}
\limsup_{n\to\infty }\frac{1}{n}\log \widetilde{\mathbb P}(X_n \in \Delta_1^c) <0 .
\end{align*}
Considering the second term in \eqref{complementprob}, we have, on $\Delta_2^c$
\[\frac{d\widetilde{\mathbb P}}{d\mathbb P} = e^{- nm + n \langle x^*, X_n\rangle    - n \mathcal L_n }\]
so that 
\[\widetilde{\mathbb P}( Y_n \in \Delta_2^c)  \leq e^{- nm + n \ell_n    - n \mathcal L_n } \,.\]
Taking the logarithm, this implies 
\begin{align*}
\limsup_{n\to \infty} \frac{1}{n}\log \widetilde{\mathbb P}( Y_n \in \Delta_2^c) & \leq -m + 
G_X(x^*) - 
G_X(x^*)-J(\varphi_m) 
\\ & =  -m - \sup_{z\in \mathcal{Y}}\{ \varphi_m(z) - \mathcal{I}_Y(z)\} \leq -m + \mathcal{I}_Y(y)
\end{align*}
which tends to $-\infty$ when $m \rightarrow \infty$.

To summarize, we have proved (\ref{complementprob}), i.e. (\ref{7.7}), which with (\ref{mino}) and (\ref{uptodelta}) gives
\[\lim_{\Delta_1 \downarrow x, \Delta_2 \downarrow y}\liminf_{n \to \infty} \frac{1}{n}\log \mathbb P(X_n \in \Delta_1, Y_n \in \Delta_2) \geq - G_X^*(x) - \mathcal I_Y (y)\,,\]
which leads to the lower bound of the LDP.
\hfill $\Box$

\section{Joint LDP for measure and truncated eigenvalues}
\label{appli}
 
In Section 4 of \cite{magicrules}, we studied the joint moment generating function of a non-negative measure $\tilde\mu^{(n)}_{I(j)}$ on a compact set $[\alpha^-,\alpha^+]$ and a collection of $j$ extremal support points $\lambda_M^\pm(j)\in \mathbb{R}^{2j}$, restricted to the compact set $[-M,M]$. For the sake of a clearer notation, we drop here the dependency on $j$. It is shown in Theorem 4.1 in \cite{magicrules}, that $\lambda_M^\pm$ satisfies the LDP with speed $n$ and good rate $\mathcal{I}_{M,\lambda^\pm}$. Furthermore, the sequence of $\tilde\mu^{(n)}_{I}$ is exponentially tight and if 
\begin{align*}
\mathcal G_n (f,s) = \mathbb E \left[\exp \left\{n \int f \, d\tilde\mu^{(n)}_{I} + n \langle s,\lambda_M^\pm\rangle \right\}\right] ,
\end{align*}
then for all $f$ such that $\log (1-f)$ is continuous and bounded and all $s\in \mathbb{R}^{2j}$, 
\begin{align} \label{oldlimit}
\lim_{n\to \infty} \frac{1}{n}\log \mathcal G_n (f,s) = G(f) + H(s) ,
\end{align}
with $G^*$ strictly convex on a set of points dense in $\{G^*<\infty\}$.

However, the rate $\mathcal{I}_{M,\lambda^\pm}$ might be non-convex and hence the dual $H^*$ is not strictly convex on a dense set. The convergence in \eqref{oldlimit} is therefore not enough to conclude the joint LDP for $(\tilde\mu^{(n)}_{I},\lambda_M^\pm)$ directly from the classical Theorem \ref{Baldi}.

To show the joint LDP, we will apply Theorem \ref{newgeneral}. Indeed,  let $D$ be the set of bounded continuous functions $f$ from $[\alpha^-,\alpha^+]$ to $\mathbb R$ such that $\sup_x f(x)<1$. If we define for $\varphi:\mathbb{R}^{2j}\to \mathbb{R}$ continuous and bounded and $f \in D$
\begin{align*}
\hat{\mathcal G}_n (f,\varphi) = \mathbb E \left[\exp \left\{n \int f \, d\tilde\mu^{(n)}_{I} + n \varphi(\lambda_M^\pm) \right\}\right] ,
\end{align*}
then the same arguments as in Section 4 of \cite{magicrules} show 
\begin{align}\label{newlimit}
\lim_{n\to \infty} \frac{1}{n}\log \hat{\mathcal G}_n (f,\varphi)  = G(f) + J(\varphi) ,
\end{align}
where
\begin{align*}
G(f) = -\int \log (1-f) \, d\mu_V 
\end{align*}
for a probability measure $\mu_V$.
Moreover, in Section 4 of \cite{magicrules} it is shown that every measure on $[\alpha^-,\alpha^+]$ with a strictly positive continuous density $h$ with respect to $\mu_V$ is an exposed point and the exposing hyperplane is the function $f=1-h^{-1}$. Since $f$ is continuous and strictly less than 1, $\gamma f\in D$ for $\gamma>1$ small enough. By the same arguments as in \cite{grz}, any measure $\mu$ with $G^*(\mu)<\infty$ may be approximated weakly by measures $\mu_n$ with such a strictly positive continuous density such that $G^*(\mu_n)$ converges to $G^*(\mu)$. This approximation is also made more precise for matrix valued measures in \cite{GaNaRomat}.  
The assumptions of Theorem \ref{newgeneral} are then satisfied, which yields the LDP for $(\tilde\mu^{(n)}_{I},\lambda_M^\pm)$ with good rate 
\begin{align*}
\mathcal{I}(\mu,\lambda) = G^*(\mu) + \mathcal{I}_{M,\lambda^\pm} . 
\end{align*}
After taking the limit $M\to \infty$, this proves the statement of Theorem 4.2 of \cite{magicrules}.

\bibliographystyle{plain}
\bibliography{bibclean}

\begin{thebibliography}{1}

\bibitem{demboz98}
A.~Dembo and O.~Zeitouni.
\newblock {\em Large Deviations Techniques and Applications}.
\newblock Springer, 1998.

\bibitem{gamb}
H.~Dette and F.~Gamboa.
\newblock Asymptotic properties of the algebraic moment range process.
\newblock {\em Acta Mathematica Hungarica}, 116:247--264, 2007.

\bibitem{GaNaRomat}
F.~Gamboa, J.~Nagel, and A.~Rouault.
\newblock Sum rules and large deviations for spectral matrix measures.
\newblock arXiv preprint arXiv:1601.08135, 2016.

\bibitem{magicrules}
F.~Gamboa, J.~Nagel, and A.~Rouault.
\newblock Sum rules via large deviations.
\newblock {\em J. Funct. Anal.}, 270(2):509 -- 559, 2016.

\bibitem{GNROPUC}
F.~Gamboa, J.~Nagel, and A.~Rouault.
\newblock Sum rules and large deviations for spectral measures on the unit
  circle.
\newblock {\em Random Matrices Theory Appl.}, 6(1), 2017.

\bibitem{grz}
F.~Gamboa, A.~Rouault, and M.~Zani.
\newblock A functional large deviation principle for quadratic forms of
  {G}aussian stationary processes.
\newblock {\em Stat. and Probab. Letters}, 43:299--308, 1999.

\end{thebibliography}

\end{document}